\documentclass[doublespacing]{elsart}

\usepackage{amssymb,amsmath}

\def\oddg{\mathop{\rm oddg}}
\def\tr{\mathop{\rm Tr}}
\def\girth{\mathop{\rm girth}}

\begin{document}

\begin{frontmatter}

\thanks[OGS]{Research partially supported by an Ontario Graduate
Scholarship}

\title{On the extreme eigenvalues of regular graphs}
\author[Kingston]{Sebastian M. Cioab\u{a}\thanksref{OGS}}
\address[Kingston]{Department of Mathematics, Queen's University, Kingston,
Ontario K7L 3N6, Canada}

\begin{center}
\emph{To the memory of Dom de Caen}
\end{center}

\begin{abstract}
In this paper, we present an elementary proof of a theorem of Serre concerning
the greatest eigenvalues of $k$-regular graphs. We also prove an analogue of
Serre's theorem regarding the least eigenvalues of $k$-regular graphs: given
$\epsilon>0$, there exist a positive constant $c=c(\epsilon,k)$ and a
nonnegative integer $g=g(\epsilon,k)$ such that for any $k$-regular graph $X$
with no odd cycles of length less than $g$, the number of eigenvalues $\mu$ of
$X$ such that $\mu \leq -(2-\epsilon)\sqrt{k-1}$ is at least $c|X|$.  This
implies a result of Winnie Li.
\end{abstract}


\end{frontmatter}

\section{Preliminaries}

Let $X$ be a graph and let $v_{0}$ be a vertex of $X$. A \emph{closed
walk} in $X$ of length $r\geq 0$ starting at $v_{0}$ is a sequence
$v_{0},v_{1},\dots,v_{r}$ of vertices of $X$ such that $v_{r}=v_{0}$
and $v_{i-1}$ is adjacent to $v_{i}$ for $1\leq i\leq r$. For $r \geq
0$, let $\Phi_{r}(X)$ denote the number of closed walks of length $r$
in $X$. A \emph{cycle} of length $r$ in $X$ is a subgraph of $X$ whose
vertices can be labeled $v_{0},\dots,v_{r}$ such that
$v_{0},\dots,v_{r}$ is a closed walk in $X$ and $v_{i}\neq v_{j}$ for
all $i,j$ with $0\leq i<j\leq r$. The \emph{girth}, denoted
$\girth(X)$, of $X$ is the length of a smallest cycle in $X$ if such a
cycle exists and $\infty$ otherwise; the \emph{oddgirth}, denoted
$\oddg(X)$, of $X$ is the length of a smallest odd cycle in $X$ if
such a cycle exists and $\infty$ otherwise. The adjacency matrix of
$X$ is the matrix $A=A(X)$ of order $|X|$, where the $(u,v)$ entry is
$1$ if the vertices $u$ and $v$ are adjacent and $0$ otherwise. It is
a well known fact that $\Phi_{r}(X)=\tr(A^{r})$, for any $r \geq
0$. The eigenvalues of $X$ are the eigenvalues of $A$. If $X$ is
$k$-regular, then it is easy to see that $k$ is an eigenvalue of $X$
with multiplicity equal to the number of components of $X$ and that
any eigenvalue $\lambda$ of $X$ satisfies $|\lambda|\leq k$. For
$l\geq 1$, we denote by $\lambda_{l}(X)$ the $l$th greatest eigenvalue
of $X$ and by $\mu_{l}(X)$ the $l$th least eigenvalue of $X$.

\section{An elementary proof of Serre's theorem}

Serre has proved the following theorem (see
\cite{DaSaVa,FeLi,LiBook,Serre1}) using Chebyschev polynomials. See
also \cite{Ci} for related results. In this section, we present an
elementary proof of Serre's result.
\begin{thm}\label{Serre}
For each $\epsilon > 0$, there exists a positive constant $c=c(\epsilon,k)$
such that for any $k$-regular graph $X$, the number of eigenvalues $\lambda$ of
$X$ with $\lambda \geq (2-\epsilon)\sqrt{k-1}$ is at least $c|X|$.
\end{thm}
For the proof of this theorem we require the next lemma which can be
deduced from McKay's work \cite{McKay}, Lemma 2.1. For the sake of
completeness, we include a short proof here.
\begin{lem}\label{McKay}
Let $v_{0}$ be a vertex of a $k$-regular graph $X$. Then the number of
closed walks of length $2s$ in $X$ starting at $v_{0}$ is greater than
or equal to $\frac{1}{s+1}{2s \choose s}k(k-1)^{s-1}$.
\end{lem}
\begin{pf*}{Proof of Lemma \ref{McKay}}
The number of closed walks of length $2s$ in $X$ starting at $v_{0}$
is at least the number of closed walks of length $2s$ starting at a
vertex $u_{0}$ in the infinite $k$-regular tree. To each closed walk
in the infinite $k$-regular tree, there corresponds a sequence of
nonnegative integers $\delta_{1},\dots, \delta_{2s}$, where
$\delta_{i}$ is the distance from $u_{0}$ after $i$ steps. The number
of such sequences is the $s$-th Catalan number
$\frac{1}{s+1}{2s\choose s}$. For each sequence of distances, there
are at least $k(k-1)^{s-1}$ closed walks of length $2s$ since for each
step away from $u_{0}$ there are $k-1$ choices ($k$ if the walk is at
$u_{0}$).\end{pf*}\qed

By Stirling's bound on $s!$ or by a simple induction argument it is easy to
see that ${2s \choose s} \geq \frac{4^{s}}{s+1}$, for any $s \geq 1$. Hence,
for any $k$-regular graph $X$ and for any $s\geq 1$, we have by Lemma
\ref{McKay}
\begin{align}\label{tr}
\tr(A^{2s})&\geq |X|\frac{1}{s+1}{2s \choose s}k(k-1)^{s-1}>
|X|\frac{1}{(s+1)^{2}}(2\sqrt{k-1})^{2s}
\end{align}
\begin{pf*}{Proof of Theorem \ref{Serre}}
Let $X$ be $k$-regular graph of order $n$ with eigenvalues $k=\lambda_{1}\geq
\dots \geq \lambda_{n} \geq -k$. Given $\epsilon > 0$, let $m$ be the number of
eigenvalues $\lambda$ of $X$ with $\lambda \geq (2-\epsilon)\sqrt{k-1}$. Then
$n-m$ of the eigenvalues of $X$ are less than $(2-\epsilon)\sqrt{k-1}$. Thus
\begin{align*}
\tr(kI+A)^{2s}&=\sum_{i=1}^{n}(k+\lambda_{i})^{2s} \\
&<(n-m)(k+(2-\epsilon)\sqrt{k-1})^{2s}+m(2k)^{2s} \\
&=m((2k)^{2s}-(k+(2-\epsilon)\sqrt{k-1})^{2s})+n(k+(2-\epsilon)\sqrt{k-1})^{2s}
\end{align*}
On the other hand, the binomial expansion and relation \eqref{tr} give
\begin{align*}
\tr(kI+A)^{2s}&=\sum_{i=0}^{2s}{2s \choose i}k^{i}\tr(A^{2s-i})\\
&\geq \sum_{j=0}^{s}{2s \choose 2j}k^{2j}\tr(A^{2s-2j})\\
&> \frac{n}{(s+1)^{2}}\sum_{j=0}^{s}{2s \choose2j}k^{2j}(2\sqrt{k-1})^{2s-2j}\\
&=\frac{n}{2(s+1)^{2}}((k+2\sqrt{k-1})^{2s}+(k-2\sqrt{k-1})^{2s})\\
&>\frac{n}{2(s+1)^{2}}(k+2\sqrt{k-1})^{2s}
\end{align*}
Thus,
\begin{equation*}
\frac{m}{n} >
\frac{\frac{1}{2(s+1)^{2}}(k+2\sqrt{k-1})^{2s}-(k+(2-\epsilon)\sqrt{k-1})^{2s}}{(2k)^{2s}-(k+(2-\epsilon)\sqrt{k-1})^{2s}}
\end{equation*}
for any $s\geq 1$. Since
\begin{align*}
\lim_{s\rightarrow \infty}
\left(\frac{(k+2\sqrt{k-1})^{2s}}{2(s+1)^{2}}\right)^{\frac{1}{2s}}&=k+2\sqrt{k-1}\\
&>k+(2-\epsilon)\sqrt{k-1}=\lim_{s\rightarrow
\infty}\left(2(k+(2-\epsilon)\sqrt{k-1})^{2s}\right)^{\frac{1}{2s}}
\end{align*}
it follows that there exists $s_{0}=s_{0}(\epsilon,k)$ such that for all $s
\geq s_{0}$
\begin{equation*}
\frac{(k+2\sqrt{k-1})^{2s}}{2(s+1)^{2}}-(k+(2-\epsilon)\sqrt{k-1})^{2s}
>(k+(2-\epsilon)\sqrt{k-1})^{2s}
\end{equation*}
Hence, if
\begin{equation*}
c(\epsilon,k)=\frac{(k+(2-\epsilon)\sqrt{k-1})^{2s_{0}}}{(2k)^{2s_{0}}-(k+(2-\epsilon)\sqrt{k-1})^{2s_{0}}}
\end{equation*}
then $c(\epsilon,k)>0$ and $m>c(\epsilon,k)n.$ \end{pf*}\qed

The proofs of Serre's theorem given in \cite{DaSaVa,FeLi,LiBook} don't
allow an easy estimation of the constant $c(\epsilon,k)$ in terms of
$\epsilon$ and $k$. We relegate the detailed analysis of the constant
obtained by those arguments to a future work \cite{CiMu}. We should
mention that Serre's theorem can be also deduced from the work of
Friedman \cite{Fri93} or Nilli \cite{Nilli2}. Friedman's results imply
an estimate of $\left(\frac{1}{2}\right)^{O\left(\frac{\log
k}{\sqrt{\epsilon}}\right)}$ for the proportion of the eigenvalues
that are at least $(2-\epsilon)\sqrt{k-1}$. Nilli's work provides a
bound of $\left(\frac{1}{2}\right)^{O\left(\frac{\log
k}{\arccos(1-\epsilon)}\right)}$. Their methods provide better bounds
on $c(\epsilon,k)$ than ours. From our proof of Serre's theorem, we
obtain that a proportion of
$\left(\frac{1}{2}\right)^{O\left(\frac{\sqrt{k}}{\epsilon}\log\left({\frac{\sqrt{k}}{\epsilon}}\right)\right)}$
of the eigenvalues are at least $(2-\epsilon)\sqrt{k-1}$. This is
because in Theorem 1 we pick $s_{0}$ such that $\frac{s_{0}}{\log
s_{0}}=\Theta\left(\frac{\sqrt{k}}{\epsilon}\right)$.

Theorem \ref{Serre} has the following consequence regarding the
asymptotics of the greatest eigenvalues of $k$-regular graphs.
\begin{cor}\label{corSerre}
Let $(X_{i})_{i\geq 0}$ be a sequence of $k$-regular graphs such that\\
$\displaystyle\lim_{i\rightarrow \infty}|X_{i}|=\infty$. Then for each
$l\geq 1$,
\begin{equation*}
\liminf_{i\rightarrow \infty}\lambda_{l}(X_{i}) \geq 2\sqrt{k-1}
\end{equation*}
\end{cor}
This corollary has also been proved directly by Serre in an appendix
to \cite{Li} using the eigenvalue distribution theorem in
\cite{Serre}. When $l=2$, we obtain the asymptotic version of the
Alon-Boppana theorem (see \cite{Alon,LPS,Nilli,Pizer} for more
details).

\section{Analogous theorems for the least eigenvalues of regular graphs}

The analogous result to Theorem \ref{Serre} for the least eigenvalues
of a $k$-regular graph is not true. For example, the eigenvalues of
line graphs are all at least $-2$. However, by adding an extra
condition to the hypothesis of Theorem \ref{Serre}, we can prove an
analogue of Serre's theorem for the least eigenvalues of a $k$-regular
graph.
\begin{thm}\label{oddg}
For any $\epsilon > 0$, there exist a positive constant
$c=c(\epsilon,k)$ and a non-negative integer $g=g(\epsilon,k)$ such
that for any $k$-regular graph $X$ with $\oddg(X)>g$, the number of
eigenvalues $\mu$ of $X$ with $\mu \leq -(2-\epsilon)\sqrt{k-1}$ is at
least $c|X|$.
\end{thm}
\begin{pf*}{Proof}
Let $X$ be a $k$-regular graph of order $n$ with eigenvalues $-k\leq
\mu_{1}\leq \mu_{2} \leq \dots \leq \mu_{n}=k$. Given $\epsilon>0$,
let $m$ be the number of eigenvalues $\mu$ of $X$ with $\mu \leq -
(2-\epsilon)\sqrt{k-1}$. Then $n-m$ of the eigenvalues of $X$ are
greater than $-(2-\epsilon)\sqrt{k-1}$. Thus
\begin{align*}
\tr(kI-A)^{2s}&=\sum_{i=1}^{n}(k-\mu_{i})^{2s}<(n-m)(k+(2-\epsilon)\sqrt{k-1})^{2s}+m(2k)^{2s}\\
&=m((2k)^{2s}-(k+(2-\epsilon)\sqrt{k-1})^{2s})+n(k+(2-\epsilon)\sqrt{k-1})^{2s}
\end{align*}
In the previous section, we proved that there exists
$s_{0}=s_{0}(\epsilon,k)$ such that for all $s \geq s_{0}$
\begin{equation*}
\frac{(k+2\sqrt{k-1})^{2s_{0}}}{2(s_{0}+1)^{2}}-(k+(2-\epsilon)\sqrt{k-1})^{2s_{0}}
>(k+(2-\epsilon)\sqrt{k-1})^{2s_{0}}
\end{equation*}
Let $g(\epsilon,k)=2s_{0}$. If $\oddg(X)>2s_{0}$, then for $0\leq j
\leq s_{0}-1$, the number of closed walks of length $2s_{0}-2j-1$ in
$X$ is $0$.  Hence, $\tr(A^{2s_{0}-2j-1})=0$, for $0\leq j\leq
s_{0}-1$. Using also \eqref{tr}, we obtain
\begin{align*}
\tr(kI-A)^{2s_{0}}&=\sum_{j=0}^{s_{0}}{2s_{0} \choose
2j}k^{2j}\tr(A^{2s_{0}-2j})-\sum_{j=0}^{s_{0}-1}{2s_{0} \choose
2j+1}k^{2j+1}\tr(A^{2s_{0}-2j-1})\\ 
&=\sum_{j=0}^{s_{0}}{2s_{0}\choose 2j}k^{2j}\tr(A^{2s_{0}-2j})>
\frac{n}{(s_{0}+1)^{2}}\sum_{j=0}^{2s_{0}}{2s_{0}\choose
2j}k^{2j}(2\sqrt{k-1})^{2s_{0}-2j}\\
&>\frac{n}{2(s_{0}+1)^{2}}(k+2\sqrt{k-1})^{2s_{0}}.
\end{align*} 
From the previous inequalities, it follows that if
\begin{equation*}
c(\epsilon,k)=\frac{(k+(2-\epsilon)\sqrt{k-1})^{2s_{0}}}{(2k)^{2s_{0}}-(k+(2-\epsilon)\sqrt{k-1})^{2s_{0}}}
\end{equation*}
then $c(\epsilon,k)>0$ and $m>c(\epsilon,k)n.$\end{pf*}\qed

The next result is an immediate consequence of Theorem \ref{oddg}.
\begin{cor}\label{coroddg}
Let $(X_{i})_{i\geq 0}$ be a sequence of $k$-regular graphs such that \\
$\lim_{i\rightarrow \infty}\oddg(X_{i})=\infty$. Then for each $l\geq 1$
\begin{equation*}
\limsup_{i\rightarrow \infty} \mu_{l}(X_{i}) \leq -2\sqrt{k-1}
\end{equation*}
\end{cor}

When $l=1$, we get the main result from \cite{Li}. Also, Corollary
\ref{coroddg} holds when $l=1$ and $\lim_{i\rightarrow
\infty}\girth(X_{i})=\infty$. This special case of Corollary
\ref{coroddg} was proved directly in \cite{LiSole} using orthogonal
polynomials and is also a consequence of the eigenvalue distribution
theorem from \cite{McKay}.

A theorem stronger than Corollary \ref{coroddg} has been proved by
Serre in \cite{Li} using the eigenvalue distribution results from
\cite{Serre}. We now present an elementary proof of this theorem.  For
$r\geq 0$, let $c_{r}(X)$ be the number of cycles of length $r$ in a
graph $X$.
\begin{thm}
Let $(X_{i})_{i\geq 0}$ be a sequence of $k$-regular graphs such
that\\ $\displaystyle\lim_{i\rightarrow\infty}|X_{i}|=\infty$. If
$\displaystyle\lim_{i\rightarrow\infty}\frac{c_{2r+1}(X_{i})}{|X_{i}|}=0$
for each $r\geq 1$, then for each $l\geq 1$
\begin{equation*}
\limsup_{i\rightarrow\infty}\mu_{l}(X_{i}) \leq -2\sqrt{k-1}
\end{equation*}
\end{thm}
\begin{pf*}{Proof}
Let $l\geq 1$. For a graph $X$ and $r\geq 1$, let $n_{2r+1}(X)$ denote
the number of vertices $v_{0}$ in the graph $X$ such that the subgraph
of $X$ induced by the vertices at distance at most $r$ from $v_{0}$ is
bipartite.  Thus, $|X|-n_{2r+1}(X)$ is the number of vertices $u_{0}$
of $X$ such that the subgraph of $X$ induced by the vertices at
distance at most $r$ from $u_{0}$ contains at least one odd
cycle. Since each such vertex is no further than $r$ from each of the
vertices of an odd cycle of length at most $2r+1$, it follows that
\begin{equation*}
|X|-n_{2r+1}(X) \leq \sum_{l=1}^{r-1}\alpha_{l,r}c_{2l+1}(X)
\end{equation*}
where $0\leq \alpha_{l,r}\leq 3(2l+1)(k-1)^{r}$. Thus, we have the
following inequalities
\begin{equation*}
1-\sum_{l=1}^{r-1}\alpha_{l,r}\frac{c_{2l+1}(X_{i})}{|X_{i}|} \leq
\frac{n_{2r+1}(X_{i})}{|X_{i}|} \leq 1
\end{equation*}
for all $r\geq 1, i\geq 0$. Hence, for each $r\geq 1$
\begin{equation}\label{n_2r}
\lim_{i\rightarrow \infty} \frac{n_{2r+1}(X_{i})}{|X_{i}|}=1
\end{equation}
For $i\geq 0$, let $A_{i}=A(X_{i})$. Then, for $i\geq 0$ and $r\geq
1$, we have
\begin{equation}\label{tr_2r}
\tr(A_{i}^{2r+1})=n_{2r+1}(X_{i})\cdot
0+(|X_{i}|-n_{2r+1}(X_{i}))\theta_{2r+1}(X_{i})
\end{equation}
where $0\leq \theta_{2r+1}(X_{i})\leq k^{2r+1}$. From \eqref{n_2r} and
\eqref{tr_2r}, we obtain that for each $r\geq 1$
\begin{equation}\label{a_2r1}
\lim_{i\rightarrow \infty} \frac{\tr(A_{i}^{2r+1})}{|X_{i}|}=0
\end{equation}
By using relation \eqref{tr}, it follows that for each $r\geq 1$
\begin{equation}\label{a_2r}
\liminf_{i\rightarrow \infty}\frac{\tr(A_{i}^{2r})}{|X_{i}|}\geq
\frac{(2\sqrt{k-1})^{2r}}{(r+1)^{2}}
\end{equation}
Now for each $i\geq 0$, we have
\begin{equation*}
\tr(kI-A_{i})^{2s}=\sum_{j=1}^{|X_{i}|}(k-\lambda_{j}(X_{i}))^{2s}\leq
(|X_{i}|-l)(k-\mu_{l}(X_{i}))^{2s}+l(2k)^{2s}
\end{equation*}
Once again, the binomial expansion gives us
\begin{equation*}
\tr(kI-A_{i})^{2s}=\sum_{j=0}^{2s} {2s \choose
j}k^{j}(-1)^{2s-j}\tr(A_{i}^{2s-j})
\end{equation*}
From the previous two relations, we get that
\begin{equation*}
(k-\mu_{l}(X_{i}))^{2s}+\frac{4^{s}lk^{2s}}{|X_{i}|-l} \geq \sum_{j=0}^{2s}{2s
\choose j}k^{j}(-1)^{2s-j}\frac{\tr(A_{i}^{2s-j})}{|X_{i}|-l}
\end{equation*}
Using relations \eqref{a_2r1} and \eqref{a_2r}, it follows that
\begin{align*}
k-\limsup_{i\rightarrow \infty} \mu_{l}(X_{i}) &\geq \left(\sum_{j=0}^{s}{2s
\choose
2j}k^{2j}\frac{(2\sqrt{k-1})^{2s-2j}}{(s-j+1)^{2}}\right)^{\frac{1}{2s}}\\
&>\left(\frac{1}{(s+1)^{2}}\sum_{j=0}^{s}{2s \choose
2j}k^{2j}(2\sqrt{k-1})^{2s-2j}\right)^{\frac{1}{2s}}\\
&>\left(\frac{1}{2(s+1)^{2}}\right)^{\frac{1}{2s}}(k+2\sqrt{k-1})
\end{align*}
for any $s\geq 1$. By taking the limit as $s\rightarrow \infty$, we get
\begin{equation*}
k-\limsup_{i\rightarrow \infty} \mu_{l}(X_{i}) \geq k+2\sqrt{k-1}
\end{equation*}
which implies the inequality stated in the theorem.\end{pf*}\qed

\section*{Acknowledgments}
This paper is part of my Ph.D. thesis at Queen's University. I am
grateful to my thesis advisors: David Gregory, Ram Murty and David
Wehlau for their help and support and to Chris Godsil and Shlomo Hoory
for their comments. I thank the referees for many useful suggestions.


\begin{thebibliography}{99}

\bibitem{Alon} N. Alon, Eigenvalues and Expanders, {\em
Combinatorica}, {\bf 6} (1986), 83-96.

\bibitem{Ci} S. M. Cioab\u{a}, Eigenvalues, Expanders and Gaps between
Primes, Ph.D. Thesis, Queen's University at Kingston, submitted
(2005).

\bibitem{CiMu} S. M. Cioab\u{a} and R. Murty, Expander Graphs and Gaps
between Primes, in progress.

\bibitem{DaSaVa} G. Davidoff, P. Sarnak and A. Vallete, {\em
Elementary Number Theory, Group Theory and Ramanujan Graphs},
Cambridge University Press, (2003).

\bibitem{FeLi} K. Feng and W.-C. Winnie Li, Spectra of hypergraphs and
applications, {\em J. of Number Theory}, {\bf 60} (1996), no.1, 1-22.

\bibitem{Fri93} J. Friedman, Some geometric aspects of graphs and
their eigenfunctions, {\em Duke Math. J.}, {\bf 69} (1993), 487-525.

\bibitem{LiBook} W.-C. Winnie Li, Number Theory with Applications,
Series of University Mathematics, Vol.7, World Scientific, (1996).

\bibitem{Li} W.-C. Winnie Li (with an appendix by J.-P. Serre), On
negative eigenvalues of regular graphs, {\em Comptes Rendus de
l'Acad$\acute{e}$}{\em mie des Sciences}, {\bf 333} (2001), issue 10,
907-912.

\bibitem{LiSole} W.-C. Winnie Li and P. Sol\'{e}, Spectra of regular
graphs and hypergraphs and orthogonal polynomials, {\em
Europ.J.Combin.}, {\bf 17} (1996), 461-477.

\bibitem{LPS} A. Lubotzky, R. Phillips and P. Sarnak, Ramanujan
Graphs, {\em Combinatorica}, {\bf 8} (1988), no.3, 261-277.

\bibitem{McKay} B. McKay, The expected eigenvalue distribution of a
large regular graph, {\em Linear Algebra and its Applications}, {\bf
40} (1981), 203-216.

\bibitem{Nilli} A. Nilli, On the second eigenvalue of a graph, {\em
Discrete Mathematics}, {\bf 91} (1991), 207-210.

\bibitem{Nilli2} A. Nilli, Tight estimates for eigenvalues of regular
graphs, {\em Electronic Journal of Combinatorics}, {\bf 11} (2004),
N9.

\bibitem{Pizer} A. Pizer, Ramanujan Graphs, {\em Computational
perspectives on number theory} (Chicago, IL, 1995), 159-178, AMS/IP
Stud. Adv.  Math., 7, Amer. Math. Soc., Providence, RI, (1998).

\bibitem{Serre1} J.-P. Serre, Private letters to W. Li dated October
8, 1990 and November 5, 1990.

\bibitem{Serre} J.-P. Serre, R\'{e}partition asymptotique des valeurs
propres de l'op\'{e}rateur de Hecke $T_{p}$, {\em J.Amer.Math.Soc.},
{\bf 10} (1997), no.1, 75-102.



\end{thebibliography}
\end{document}